\documentclass[12pt]{article}%
\usepackage{amsmath}
\usepackage{amsfonts}
\usepackage{amssymb}
\usepackage{graphicx}%
\providecommand{\U}[1]{\protect\rule{.1in}{.1in}}
\newtheorem{theorem}{Theorem}[section]

\numberwithin{equation}{section}

\begin{document}

\title{Implementation of the inverse scattering transform method for the nonlinear
Schr\"{o}dinger equation}
\author{Vladislav V. Kravchenko\\{\small Departamento de Matem\'{a}ticas, Cinvestav, Unidad Quer\'{e}taro, }\\{\small Libramiento Norponiente \#2000, Fracc. Real de Juriquilla,
Quer\'{e}taro, Qro., 76230 MEXICO.}\\{\small e-mail: vkravchenko@math.cinvestav.edu.mx}}
\maketitle

\begin{abstract}
We study the initial-value problem for the nonlinear Schr\"{o}dinger equation.
Application of the inverse scattering transform method involves solving direct
and inverse scattering problems for the Zakharov-Shabat system with complex
potentials. We solve these problems by using new series representations for
the Jost solutions of the Zakharov-Shabat system. The representations have the
form of power series with respect to a transformed spectral parameter.\ In
terms of the representations, solution of the direct scattering problem
reduces to computing the series coefficients following a simple recurrent
integration procedure, computation of the scattering coefficients by
multiplying corresponding pairs of polynomials (partial sums of the series
representations) and locating zeros of a polynomial inside the unit disk.

Solution of the inverse scattering problem reduces to the solution of a system
of linear algebraic equations for the power series coefficients, while the
potential is recovered from the first coefficients. The system is obtained
directly from the scattering relations. Thus, unlike other existing
techniques, the method does not involve solving the Gelfand-Levitan-Marchenko
equation or the matrix Riemann-Hilbert problem.

The overall approach leads to a simple and efficient algorithm for the
numerical solution of the initial-value problem for the nonlinear
Schr\"{o}dinger equation, which is illustrated by numerical examples.\ 

\end{abstract}

\section{Introduction}

We consider the inverse scattering transform method (ISTM) for solving the
initial-value problem for the nonlinear Schr\"{o}dinger equation (NLSE)
\begin{equation}
iq_{t}+q_{xx}+2q\left\vert q\right\vert ^{2}=0, \label{NLSE intro}%
\end{equation}%
\begin{equation}
q(x,0)=q_{0}(x). \label{initial cond}%
\end{equation}

We assume that $q_{0}$ belongs to the class $\mathcal{Q}$ of\ complex valued
functions, such that for some $k\geq1$ the expression $(1+\left\vert
x\right\vert ^{k})\left\vert q_{0}(x)\right\vert $ is both integrable and
square integrable over $-\infty<x<\infty$. In this case the ISTM is applicable
to (\ref{NLSE intro}), (\ref{initial cond}), see \cite{Zhou 1989} and
discussion in \cite[subsect. 5.1.6]{Shaw}. Application of ISTM to
(\ref{NLSE intro}), (\ref{initial cond}) involves solving a direct scattering
problem for the Zakharov-Shabat system
\begin{align}
\frac{dn_{1}(x)}{dx}+i\rho n_{1}(x)  &  =q(x)n_{2}(x),\label{ZS1}\\
\frac{dn_{2}(x)}{dx}-i\rho n_{2}(x)  &  =-\overline{q}(x)n_{1}(x),\quad
x\in(-\infty,\infty), \label{ZS2}%
\end{align}
with the potential $q=q_{0}$ to obtain the initial scattering data, and
solving the inverse scattering problem for the Zakharov-Shabat system, thus
reconstructing $q(x,t)$ from the evolved scattering data. Here $\rho
\in\mathbb{C}$ is a spectral parameter. The procedure resembles the direct and
inverse Fourier transform and therefore, it is often referred to as the
nonlinear Fourier transform. Its theory is presented in a number of books and
review papers (see, e.g., \cite{Ablowitz Segur}, \cite{Gerdjikov et al 2008},
\cite{Lamb}, \cite{Yousefi et al 2014 IEEE}). Numerical techniques are
discussed in dozens of papers (see \cite{Arico et al 2011}, \cite{Civelli et
al 2015}, \cite{Delitsyn 2022}, \cite{Fermo et al 2016}, \cite{Frangos et al
1991}, \cite{Frumin et al 2015}, \cite{Gorbenko et al}, \cite{Le et al 2014},
\cite{Medvedev et al 2023}, \cite{Mullyadzhanov et al 2021}, \cite{Tang et
al}, \cite{Trogdon2021}, \cite{Trogdon Olver}, \cite{Turitsyn et al 2017},
\cite{Xiao et al 2002}). While the direct scattering is realized numerically
by applying the factorization property of the transfer (scattering) matrix
and/or discretization of the Zakharov-Shabat system, the numerical inverse
scattering is based on the solution of the Gelfand-Levitan integral equation
or on the solution of a related matrix Riemann-Hilbert problem. A review of
the existing numerical techniques can be consulted in \cite{Delitsyn 2022},
\cite{Medvedev et al 2023}, \cite{Trogdon2021}, \cite{Turitsyn et al 2017}.

In this work we develop an entirely different approach initiated in
\cite{Kr2025ZSreal}, \cite{Kr2025AKNS}. Its main feature consists in using
certain series representations for the Jost solutions of the Zakharov-Shabat
system. The series are power series with respect to the transformed spectral
parameter
\begin{equation}
z=\frac{\frac{1}{2}+i\rho}{\frac{1}{2}-i\rho}, \label{z Intro}%
\end{equation}
and therefore are called spectral parameter power series (SPPS). For computing
the coefficients of the series there is a simple recurrent integration
procedure \cite{Kr2025AKNS}, which we adapt here to the solutions of
(\ref{ZS1}), (\ref{ZS2}). The procedure allows one to compute hundreeds of the
coefficients and can be easily implemented using any available numerical
computing environment.

Since (\ref{z Intro}) is a M\"{o}bius transformation of the upper half-plane
(where $\rho$ is considered) to the unit disk, all the computations of the
scattering data in terms of $z$ are reduced to the unit disk, including the
location of the eigenvalues. Moreover, since we deal with the power series in
$z$, in practice the computation of the scattering data reduces to
manipulations with polynomials in $z$, which are partial sums of the power
series representations. Numerical tests reveal extreme accuracy and speed of
this approach.

The method of solution of the inverse scattering problem with the aid of the
SPPS representations is also elementary. We substitute the SPPS for the Jost
solutions into the scattering relations and write down a resulting system of
linear algebraic equations for the SPPS coefficients. Moreover, we show that
the potential can be reconstructed from the very first coefficient, so there
is no necessity in computing many of them. This observation greatly reduces
the number of the unknowns\ in the system of linear algebraic equations, thus
simplifying the solution of the inverse problem.

The method is simple, fast and accurate, which is illustrated by several
numerical examples.

Besides this introduction the paper contains seven sections. Section
\ref{Sect SPPS ZS} contains some definitions and notations related to the
Zakharov-Shabat system. In Section \ref{Sect NLSE by ISTM} we briefly recall
the inverse scattering transform method. In Section \ref{Sect SPPS ZS} we
introduce the SPPS representations for the Jost solutions together with the
recurrent integration procedure for the coefficients of the series. In
Sections \ref{SectDirect} and \ref{Sect Inverse} we introduce the methods for
solving the direct and inverse scattering problems, respectively. Section
\ref{Sect Numerics} presents some numerical experiments.

\section{Zakharov-Shabat system, Jost solutions and scattering
data\label{Sect ZS Jost SD}}

Let $q\in\mathcal{Q}$. Consider the Zakharov-Shabat system (\ref{ZS1}),
(\ref{ZS2}). It is convenient to regard its solutions as functions of $\rho$
as well, so we write
\[
n=n(\rho,x)=\left(
\begin{tabular}
[c]{c}%
$n_{1}(\rho,x)$\\
$n_{2}(\rho,x)$%
\end{tabular}
\right)  .
\]

System (\ref{ZS1}), (\ref{ZS2}) possesses so-called Jost solutions, satisfying
prescribed conditions at $\pm\infty$, respectively:
\[
\varphi(\rho,x)\sim\left(
\begin{tabular}
[c]{c}%
$1$\\
$0$%
\end{tabular}
\right)  e^{-i\rho x},\quad\widetilde{\varphi}(\rho,x)=\left(
\begin{tabular}
[c]{c}%
$0$\\
$-1$%
\end{tabular}
\right)  e^{i\rho x},\quad x\rightarrow-\infty
\]
and%
\[
\psi(\rho,x)\sim\left(
\begin{tabular}
[c]{c}%
$0$\\
$1$%
\end{tabular}
\right)  e^{i\rho x},\quad\widetilde{\psi}(\rho,x)\sim\left(
\begin{tabular}
[c]{c}%
$1$\\
$0$%
\end{tabular}
\right)  e^{-i\rho x},\quad x\rightarrow\infty.
\]
The relations hold%
\begin{equation}
\widetilde{\varphi}(\rho,x)=\left(
\begin{tabular}
[c]{c}%
$\overline{\varphi}_{2}(\overline{\rho},x)$\\
$-\overline{\varphi}_{1}(\overline{\rho},x)$%
\end{tabular}
\right)  ,\quad\widetilde{\psi}(\rho,x)=\left(
\begin{tabular}
[c]{c}%
$\overline{\psi}_{2}(\overline{\rho},x)$\\
$-\overline{\psi}_{1}(\overline{\rho},x)$%
\end{tabular}
\right)  . \label{phitil and psitil}%
\end{equation}
The pairs of the solutions $\varphi(\rho,x)$, $\widetilde{\varphi}(\rho,x)$
and $\psi(\rho,x)$, $\widetilde{\psi}(\rho,x)$ are linearly independent, and
hence there exist the scalars $\mathbf{a}(\rho)$ and $\mathbf{b}(\rho)$,
called scattering coefficients, such that for all $\rho\in\mathbb{R}$ we have
the scattering relations
\begin{equation}
\varphi(\rho,x)=\mathbf{a}(\rho)\widetilde{\psi}(\rho,x)+\mathbf{b}(\rho
)\psi(\rho,x) \label{rel 1}%
\end{equation}
and
\begin{equation}
\widetilde{\varphi}(\rho,x)=\overline{\mathbf{b}}(\rho)\widetilde{\psi}%
(\rho,x)-\overline{\mathbf{a}}(\rho)\psi(\rho,x). \label{rel 2}%
\end{equation}

The identity is valid%
\begin{equation}
\mathbf{a}(\rho)\overline{\mathbf{a}}(\rho)+\mathbf{b}(\rho)\overline
{\mathbf{b}}(\rho)=1. \label{aatil}%
\end{equation}

In terms of the Jost solutions the scattering coefficients can be written as
follows \cite{Ablowitz Segur}%
\begin{equation}
\mathbf{a}(\rho)=W\left[  \varphi(\rho,0);\psi(\rho,0)\right]  ,\qquad
\mathbf{b}(\rho)=-W\left[  \varphi(\rho,0);\widetilde{\psi}(\rho,0)\right]  ,
\label{atila}%
\end{equation}
where $W$ denotes the Wronskian, $W\left[  \varphi;\psi\right]  :=\varphi
_{1}\psi_{2}-\varphi_{2}\psi_{1}$.

System (\ref{ZS1}), (\ref{ZS2}) may admit a number of the eigenvalues, which
coincide with zeros of $\mathbf{a}(\rho)$ in $\mathbb{C}^{+}:=\left\{  \rho
\in\mathbb{C}\mid\operatorname{Im}\rho>0\right\}  $. We will assume that the
number of the eigenvalues is finite. If $\mathbf{a}(\rho_{m})=0$,
$\operatorname{Im}\rho_{m}>0$, then the corresponding Jost solutions
$\varphi(\rho_{m},x)$ and $\psi(\rho_{m},x)$ are linearly dependent, hence
there exists a multiplier constant $c_{m}\neq0$ such that
\begin{equation}
\varphi(\rho_{m},x)=c_{m}\psi(\rho_{m},x). \label{eig1}%
\end{equation}
$c_{m}$ is called the norming constant associated with the eigenvalue
$\rho_{m}$. The number of the eigenvalues we denote by $M$.

The direct scattering problem consists in finding $\mathbf{a}(\rho)$,
$\mathbf{b}(\rho)$ for all $\rho\in\mathbb{R}$ as well as in finding all the
eigenvalues and norming constants. In other words, when solving the direct
scattering problem for system (\ref{ZS1}), (\ref{ZS2}) one needs to find the
set of the scattering data%
\begin{equation}
SD:=\left\{  \mathbf{a}(\rho),\,\mathbf{b}(\rho),\,\rho\in\mathbb{R};\left\{
\rho_{m},c_{m}\right\}  _{m=1}^{M}\right\}  . \label{SD}%
\end{equation}

The inverse scattering problem consists in recovering the potential $q(x)$ in
(\ref{ZS1}), (\ref{ZS2}) for all $x\in\mathbb{R}$ from the set of the
scattering data $SD$.

\section{Solution of the NLSE by the ISTM, evolution of scattering
data\label{Sect NLSE by ISTM}}

To find a solution $q(x,t)$ of the initial-value problem (\ref{NLSE intro}),
(\ref{initial cond}) one can solve a direct scattering problem for the
potential $q_{0}(x)$, which gives a set of the scattering data
\[
SD_{0}=\left\{  \mathbf{a}(\rho),\,\mathbf{b}(\rho),\,\rho\in\mathbb{R}%
;\,\left\{  \rho_{m},\,c_{m}\right\}  _{m=1}^{M}\right\}
\]
corresponding to the initial time $t=0$. Next, construct the scattering data
corresponding to the time $t$ as follows
\[
SD_{t}=\left\{  \mathbf{a}(\rho),\,\mathbf{b}(\rho)e^{4i\rho^{2}t},\,\rho
\in\mathbb{R};\,\left\{  \rho_{m},\,c_{m}e^{4i\rho_{m}^{2}t}\right\}
_{m=1}^{M}\right\}  .
\]
Finally, the solution of the inverse scattering problem for this set of the
scattering data is precisely $q(x,t)$. This three-step procedure is the ISTM
\cite{ZS}, \cite{Ablowitz Segur}, \cite{Shaw}, \cite{Lamb}. Unlike usual
time-stepping numerical techniques, the ISTM does not require computing the
solution of (\ref{NLSE intro}), (\ref{initial cond}) for the previous times,
but allows instead to obtain it directly for any chosen time $t$, which can be
very large.

\section{Spectral parameter power series for Jost
solutions\label{Sect SPPS ZS}}

Let us introduce the notation
\[
z=z\left(  \rho\right)  :=\frac{\frac{1}{2}+i\rho}{\frac{1}{2}-i\rho}%
,\quad\widetilde{z}=\widetilde{z}\left(  \rho\right)  :=\frac{\frac{1}%
{2}-i\rho}{\frac{1}{2}+i\rho}.
\]
For all $\rho\in\overline{\mathbb{C}^{+}}$, $z\left(  \rho\right)
\in\overline{D}=\left\{  z\in\mathbb{C}:\,\left\vert z\right\vert
\leq1\right\}  $, while $\widetilde{z}\left(  \rho\right)  \in\overline
{\widetilde{D}}=\left\{  \widetilde{z}\in\mathbb{C}:\,\left\vert \widetilde
{z}\right\vert \leq1\right\}  $.

\begin{theorem}
\cite{Kr2025AKNS} Let $q(x)$ belong to the class $\mathcal{Q}$. Then the
corresponding Jost solutions of (\ref{ZS1}), (\ref{ZS2}) admit the series
representations%
\begin{equation}
\varphi(\rho,x)=e^{-i\rho x}\left(  \left(
\begin{tabular}
[c]{c}%
$1$\\
$0$%
\end{tabular}
\ \right)  +\left(  z+1\right)  \sum_{n=0}^{\infty}\left(  -1\right)
^{n}z^{n}b_{n}(x)\right)  , \label{phi=}%
\end{equation}%
\begin{equation}
\widetilde{\varphi}(\rho,x)=e^{i\rho x}\left(  \left(
\begin{tabular}
[c]{c}%
$0$\\
$-1$%
\end{tabular}
\ \right)  +\left(  \widetilde{z}+1\right)  \sum_{n=0}^{\infty}\left(
-1\right)  ^{n}\widetilde{z}^{n}\widetilde{b}_{n}(x)\right)  , \label{phitil=}%
\end{equation}%
\begin{equation}
\psi(\rho,x)=e^{i\rho x}\left(  \left(
\begin{tabular}
[c]{c}%
$0$\\
$1$%
\end{tabular}
\ \right)  +\left(  z+1\right)  \sum_{n=0}^{\infty}\left(  -1\right)
^{n}z^{n}a_{n}(x)\right)  , \label{psi=}%
\end{equation}%
\begin{equation}
\widetilde{\psi}(\rho,x)=e^{-i\rho x}\left(  \left(
\begin{tabular}
[c]{c}%
$1$\\
$0$%
\end{tabular}
\ \right)  +\left(  \widetilde{z}+1\right)  \sum_{n=0}^{\infty}\left(
-1\right)  ^{n}\widetilde{z}^{n}\widetilde{a}_{n}(x)\right)  , \label{psitil=}%
\end{equation}
where the coefficients $a_{n}$, $\widetilde{a}_{n}$, $b_{n}$ and
$\widetilde{b}_{n}$ are vector functions, e.g.,
\[
b_{n}(x)=\left(
\begin{tabular}
[c]{c}%
$b_{1,n}(x)$\\
$b_{2,n}(x)$%
\end{tabular}
\right)  .
\]
For all $x\in\mathbb{R}$, the power series in (\ref{phi=}), (\ref{psi=})
converge in the unit disk $D=\left\{  z\in\mathbb{C}:\,\left\vert z\right\vert
<1\right\}  $ of the complex plane of the variable $z$, and the power series
in (\ref{phitil=}), (\ref{psitil=}) converge in the unit disk $\widetilde
{D}=\left\{  \widetilde{z}\in\mathbb{C}:\,\left\vert \widetilde{z}\right\vert
<1\right\}  $ of the complex plane of the variable $\widetilde{z}$.
\end{theorem}

The series representations (\ref{phi=})-(\ref{psitil=}) of the Jost solutions
are called \textbf{spectral parameter power series} (SPPS), in analogy with
those introduced in \cite{KrPorter2010} for regular solutions of
Sturm-Liouville equations.

From the relations (\ref{phitil and psitil}) it is easy to derive the
relations for the SPPS coefficients:%
\begin{equation}
\widetilde{a}_{1,n}(x)=\overline{a}_{2,n}(x),\quad\widetilde{a}_{2,n}%
(x)=-\overline{a}_{1,n}(x), \label{a1ntil=a2nbar}%
\end{equation}%
\begin{equation}
\widetilde{b}_{1,n}(x)=\overline{b}_{2,n}(x),\quad\widetilde{b}_{2,n}%
(x)=-\overline{b}_{1,n}(x). \label{b1ntil=b2nbar}%
\end{equation}
Indeed, considering, e.g., $\widetilde{\varphi}_{1}(\rho,x)$ and taking into
account that $\widetilde{z}\left(  \rho\right)  =\overline{z}\left(
\overline{\rho}\right)  $, we have
\begin{align*}
\widetilde{\varphi}_{1}(\rho,x)  &  =\overline{\varphi}_{2}(\overline{\rho
},x)=\overline{e^{-i\overline{\rho}x}}\left(  \overline{z}\left(
\overline{\rho}\right)  +1\right)  \sum_{n=0}^{\infty}\left(  -1\right)
^{n}\overline{z}^{n}\left(  \overline{\rho}\right)  \overline{b}_{2,n}(x)\\
&  =e^{i\rho x}\left(  \widetilde{z}\left(  \rho\right)  +1\right)  \sum
_{n=0}^{\infty}\left(  -1\right)  ^{n}\widetilde{z}^{n}\left(  \rho\right)
\overline{b}_{2,n}(x).
\end{align*}
Comparison of this series with (\ref{phitil=}) gives us the first equality in
(\ref{b1ntil=b2nbar}). The other equalities are proved similarly.

In \cite{Kr2025AKNS} a recurrent integration procedure was deduced for
computing the SPPS coefficients. Here we adapt it to the case of the
Zakharov-Shabat system (\ref{ZS1}), (\ref{ZS2}).

First, we note that since $z\left(  \frac{i}{2}\right)  =0$, from the series
representations (\ref{phi=}) and (\ref{psi=}) we obtain%
\begin{equation}
\varphi(\frac{i}{2},x)=e^{\frac{x}{2}}\left(  \left(
\begin{tabular}
[c]{c}%
$1$\\
$0$%
\end{tabular}
\ \ \right)  +b_{0}(x)\right)  ,\quad\psi(\frac{i}{2},x)=e^{-\frac{x}{2}%
}\left(  \left(
\begin{tabular}
[c]{c}%
$0$\\
$1$%
\end{tabular}
\ \ \right)  +a_{0}(x)\right)  . \label{phi(i/2)}%
\end{equation}

Thus, to obtain the first SPPS coefficients, it is sufficient to compute the
Jost solutions $\varphi$ and $\psi$ corresponding to $\rho=\frac{i}{2}$. Then%
\begin{equation}
a_{0}(x)=e^{\frac{x}{2}}\psi(\frac{i}{2},x)-\left(
\begin{tabular}
[c]{c}%
$0$\\
$1$%
\end{tabular}
\ \right)  , \label{a0}%
\end{equation}%
\begin{equation}
b_{0}(x)=e^{-\frac{x}{2}}\varphi(\frac{i}{2},x)-\left(
\begin{tabular}
[c]{c}%
$1$\\
$0$%
\end{tabular}
\ \right)  . \label{b0}%
\end{equation}
Denote
\[
f(x):=a_{2,0}(x)+1=e^{\frac{x}{2}}\psi_{2}(\frac{i}{2},x).
\]
Then, using (\ref{ZS1}), (\ref{ZS2}) we have
\[
f^{\prime}(x)=-e^{\frac{x}{2}}\overline{q}(x)\psi_{1}(\frac{i}{2},x).
\]
Similarly, denote
\[
g(x):=b_{1,0}(x)+1=e^{-\frac{x}{2}}\varphi_{1}(\frac{i}{2},x).
\]
Then
\[
g^{\prime}(x)=e^{-\frac{x}{2}}q(x)\varphi_{2}(\frac{i}{2},x).
\]
The subsequent SPPS coefficients $\left\{  a_{n}(x)\right\}  _{n=1}^{\infty}$
are obtained as follows. Denote%
\[
h_{n}(x):=e^{-x}\left(  a_{1,n-1}^{\prime}(x)+a_{1,n-1}(x)-q(x)a_{2,n-1}%
(x)\right)  .
\]
and
\[
H_{n}(x):=\int_{x}^{\infty}f(s)h_{n}(s)ds.
\]
Now,
\begin{equation}
a_{2,n}(x)=-f(x)\int_{x}^{\infty}\frac{e^{t}\overline{q}(t)}{f^{2}(t)}%
H_{n}(t)dt,\quad n=1,2,\ldots, \label{a1}%
\end{equation}%
\begin{equation}
a_{1,n}(x)=-\frac{f^{\prime}(x)}{\overline{q}(x)f(x)}a_{2,n}(x)-\frac{e^{x}%
}{f(x)}H_{n}(x),\quad n=1,2,\ldots\label{a3}%
\end{equation}
with
\[
a_{1,0}(x)=-\frac{f^{\prime}(x)}{\overline{q}(x)}.
\]
The derivatives $a_{1,n-1}^{\prime}$ required for calculating $H_{n}(x)$ are
computed without differentiating:
\[
a_{1,0}^{\prime}(x)=a_{1,0}(x)+q(x)(a_{2,0}(x)+1),
\]%
\begin{equation}
a_{1,n}^{\prime}(x)=a_{1,n-1}^{\prime}(x)+a_{1,n}(x)+a_{1,n-1}(x)+q(x)(a_{2,n}%
(x)-a_{2,n-1}(x)),\quad n=1,2,\ldots. \label{a4}%
\end{equation}

Similarly, the coefficients $\left\{  b_{n}(x)\right\}  _{n=1}^{\infty}$ are
computed as follows. Denote
\[
p_{n}(x):=e^{x}\left(  b_{2,n-1}^{\prime}(x)-b_{2,n-1}(x)+\overline
{q}(x)b_{1,n-1}(x)\right)
\]
and
\[
P_{n}(x):=\int_{-\infty}^{x}g(s)p_{n}(s)ds.
\]
We have
\begin{equation}
b_{1,n}(x)=g(x)\int_{-\infty}^{x}\frac{e^{-t}q(t)}{g^{2}(t)}P_{n}(t)dt,
\label{b1}%
\end{equation}%
\begin{equation}
b_{2,n}(x)=\frac{g^{\prime}(x)}{q(x)g(x)}b_{1,n}(x)+\frac{e^{-x}}{g(x)}%
P_{n}(x) \label{b3}%
\end{equation}
with
\[
b_{2,0}(x)=\frac{g^{\prime}(x)}{q(x)},
\]
and the derivatives $b_{2,n}^{\prime}$ are computed by%
\[
b_{2,0}^{\prime}(x)=-b_{2,0}(x)-\overline{q}(x)\left(  b_{1,0}(x)+1\right)  ,
\]%
\begin{equation}
b_{2,n}^{\prime}(x)=b_{2,n-1}^{\prime}(x)-b_{2,n}(x)-b_{2,n-1}(x)-\overline
{q}(x)\left(  b_{1,n}(x)-b_{1,n-1}(x)\right)  . \label{b4}%
\end{equation}

The following asymptotic behaviour of the coefficients is valid
\cite{Kr2025AKNS}. We have
\[
a_{1,n}(x)=o(1),\quad a_{1,n}^{\prime}(x)=o(1),\quad n=0,1,\ldots,
\]%
\[
a_{2,0}(x)=o(1),\quad a_{2,n}(x)=o(\int_{x}^{\infty}\overline{q}(t)dt),\quad
n=1,2,\ldots,
\]
when $x\rightarrow+\infty$, and%
\[
b_{2,n}(x)=o(1),\quad b_{2,n}^{\prime}(x)=o(1),\quad n=0,1,\ldots,
\]
\[
b_{1,0}(x)=o(1),\quad b_{1,n}(x)=o(\int_{x}^{\infty}q(t)dt),\quad
n=1,2,\ldots,
\]
when $x\rightarrow-\infty$.

The SPPS coefficients $\left\{  \widetilde{a}_{n}(x),\,\widetilde{b}%
_{n}(x)\right\}  _{n=0}^{\infty}$ are obtained from (\ref{a1ntil=a2nbar}) and
(\ref{b1ntil=b2nbar}).

\section{Solution of the direct problem\label{SectDirect}}

Substitution of the SPPS (\ref{phi=})-(\ref{psitil=}) into the first equality
in (\ref{atila}) leads to the SPPS representation for $\mathbf{a}(\rho)$:%
\begin{align}
\mathbf{a}(\rho)  &  =\varphi_{1}(\rho,0)\psi_{2}(\rho,0)-\varphi_{2}%
(\rho,0)\psi_{1}(\rho,0)\nonumber\\
&  =\left(  1+\left(  z+1\right)  \sum_{n=0}^{\infty}\left(  -1\right)
^{n}z^{n}b_{1,n}(0)\right)  \left(  1+\left(  z+1\right)  \sum_{n=0}^{\infty
}\left(  -1\right)  ^{n}z^{n}a_{2,n}(0)\right) \nonumber\\
&  -\left(  z+1\right)  ^{2}\left(  \sum_{n=0}^{\infty}\left(  -1\right)
^{n}z^{n}b_{2,n}(0)\right)  \left(  \sum_{n=0}^{\infty}\left(  -1\right)
^{n}z^{n}a_{1,n}(0)\right)  \label{SPPS a}%
\end{align}
for all $\rho\in\overline{\mathbb{C}^{+}}$

Analogously, the SPPS representation for $\mathbf{b}(\rho)$ is obtained by
substituting (\ref{phi=})-(\ref{psitil=}) into the second equality in
(\ref{atila}) and taking into account (\ref{a1ntil=a2nbar}) together with the
equality $\widetilde{z}\left(  \rho\right)  =\overline{z}\left(  \rho\right)
$ for $\rho\in\mathbb{R}$. We have
\begin{align}
\mathbf{b}(\rho)  &  =\varphi_{2}(\rho,0)\widetilde{\psi}_{1}(\rho
,0)-\varphi_{1}(\rho,0)\widetilde{\psi}_{2}(\rho,0)\nonumber\\
&  =\left(  z+1\right)  \left(  \sum_{n=0}^{\infty}\left(  -1\right)
^{n}z^{n}b_{2,n}(0)\right)  \left(  1+\left(  \overline{z}+1\right)
\sum_{n=0}^{\infty}\left(  -1\right)  ^{n}\overline{z}^{n}\overline{a}%
_{2,n}(0)\right) \nonumber\\
&  +\left(  \overline{z}+1\right)  \left(  1+\left(  z+1\right)  \sum
_{n=0}^{\infty}\left(  -1\right)  ^{n}z^{n}b_{1,n}(0)\right)  \left(
\sum_{n=0}^{\infty}\left(  -1\right)  ^{n}\overline{z}^{n}\overline{a}%
_{1,n}(0)\right)  \label{SPPS b}%
\end{align}
for all $\rho\in\mathbb{R}$.

Truncating the series in (\ref{SPPS a}) and (\ref{SPPS b}) up to a
sufficiently large number $N$ gives us a simple way of computing the
scattering coefficients. Note that when $\rho\in\mathbb{R}$, $z$ belongs to
the unit circle centered at the origin; when $\rho\in\mathbb{C}^{+}$, $z$
belongs to the unit disk $D$. This observation implies that the computation of
the eigenvalues reduces to the location of zeros of (\ref{SPPS a}) inside $D$.
Numerically, the truncated series (\ref{SPPS a}) is considered, which is
nothing but a polynomial in the variable $z$. Thus, the computation of the
eigenvalues reduces to the location of the roots of the polynomial inside the
unit disk.

The norming constants $c_{m}$ are also computed by definition:%
\[
c_{m}=\frac{\varphi_{1}(\rho_{m},0)}{\psi_{1}(\rho_{m},0)}=\frac{1+\left(
z_{m}+1\right)  \sum_{n=0}^{\infty}\left(  -1\right)  ^{n}z_{m}^{n}b_{1,n}%
(0)}{\left(  z_{m}+1\right)  \sum_{n=0}^{\infty}\left(  -1\right)  ^{n}%
z_{m}^{n}a_{1,n}(0)}.
\]
Equally, the quotients of the second components can be used. For the numerical
implementation, here again, the truncated sums should be considered.

Summarizing, the direct scattering problem is solved as follows. Given $q(x)$,
compute a sufficiently large set of the coefficients $\left\{  a_{n}%
(0),b_{n}(0)\right\}  _{n=0}^{N}$ following the recurrent integration
procedure from Section \ref{Sect SPPS ZS}. Next, the scattering coefficients
are computed for $\rho\in\mathbb{R}$ by (\ref{SPPS a}), (\ref{SPPS b}), where
the sums truncated up to $N$ are considered. Next, the eigenvalues $\rho
_{m}\in\mathbb{C}^{+}$ are computed as $\rho_{m}=\frac{z_{m}-1}{2i\left(
z_{m}+1\right)  }$, where $z_{m}$ are roots of the polynomial
\begin{align*}
\mathbf{a}_{N}(z)  &  =\left(  1+\left(  z+1\right)  \sum_{n=0}^{N}\left(
-1\right)  ^{n}z^{n}b_{1,n}(0)\right)  \left(  1+\left(  z+1\right)
\sum_{n=0}^{N}\left(  -1\right)  ^{n}z^{n}a_{2,n}(0)\right) \\
&  -\left(  z+1\right)  ^{2}\left(  \sum_{n=0}^{N}\left(  -1\right)  ^{n}%
z^{n}b_{2,n}(0)\right)  \left(  \sum_{n=0}^{N}\left(  -1\right)  ^{n}%
z^{n}a_{1,n}(0)\right)  ,
\end{align*}
which are located in $D$.

\section{Inverse scattering\label{Sect Inverse}}

First, let us notice that for recovering $q(x)$ it is sufficient to find any
of the first SPPS coefficients. For example, assume that $b_{0}(x)$ is known.
This means that the Jost solution $\varphi(\frac{i}{2},x)$ corresponding to
$\rho=\frac{i}{2}$ is known (formula (\ref{phi(i/2)})):
\[
\varphi_{1}(\frac{i}{2},x)=e^{\frac{x}{2}}\left(  1+b_{1,0}(x)\right)
,\quad\varphi_{2}(\frac{i}{2},x)=e^{\frac{x}{2}}b_{2,0}(x).
\]
From (\ref{ZS1}), (\ref{ZS2}) for $\rho=\frac{i}{2}$ we have that%
\begin{equation}
q(x)=\frac{\varphi_{1}^{\prime}(\frac{i}{2},x)-\varphi_{1}(\frac{i}{2}%
,x)/2}{\varphi_{2}(\frac{i}{2},x)},\quad\overline{q}(x)=-\frac{\varphi
_{2}^{\prime}(\frac{i}{2},x)+\varphi_{2}(\frac{i}{2},x)/2}{\varphi_{1}%
(\frac{i}{2},x)}. \label{q and r}%
\end{equation}
Thus, obtaining $b_{0}(x)$ is sufficient for recovering $q(x)$. Analogously,
$q(x)$ can be recovered from any of the first SPPS coefficients.

Next, let us construct a system of linear algebraic equations for the SPPS
coefficients. Consider the first equalities in (\ref{rel 1}) and
(\ref{rel 2}). We have%

\begin{equation}
\varphi_{1}(\rho,x)=\mathbf{a}(\rho)\overline{\psi}_{2}(\rho,x)+\mathbf{b}%
(\rho)\psi_{1}(\rho,x) \label{eq1}%
\end{equation}
and
\begin{equation}
\overline{\varphi}_{2}(\rho,x)=\overline{\mathbf{b}}(\rho)\overline{\psi}%
_{2}(\rho,x)-\overline{\mathbf{a}}(\rho)\psi_{1}(\rho,x). \label{eq2}%
\end{equation}
Note that considering the second equalities in (\ref{rel 1}) and (\ref{rel 2})
one arrives at the same pair of equations up to a complex conjugation.

Observe that in terms of the SPPS coefficients, for any $\rho\in\mathbb{R}$
and for any $x\in\mathbb{R}$ equalities (\ref{eq1}) and (\ref{eq2}) are
equations for $\left\{  a_{1,n}(x),b_{1,n}(x),\overline{a}_{2,n}%
(x),\overline{b}_{2,n}(x)\right\}  _{n=0}^{\infty}$. Additionally, if the
discrete spectrum is not empty, the equalities (\ref{eig1}) for all $\rho_{m}%
$, $m=1,\ldots,M$ can be written in the form
\begin{equation}
\varphi_{1}(\rho_{m},x)=c_{m}\psi_{1}(\rho_{m},x), \label{eq3}%
\end{equation}%
\begin{equation}
\overline{\varphi}_{2}(\rho_{m},x)=\overline{c}_{m}\overline{\psi}_{2}%
(\rho_{m},x), \label{eq4}%
\end{equation}
where (\ref{eq4}) is obtained by applying complex conjugation to the second
equation in (\ref{eig1}). Equalities (\ref{eq3}) and (\ref{eq4}) are equations
for the same four sets of the SPPS coefficients.

Let us consider equalities (\ref{eq1}) and (\ref{eq2}) at a sufficiently large
number $K$ of points $\rho_{k}\in\mathbb{R}$. Substituting the truncated SPPS
representations for $\varphi_{1}$, $\overline{\varphi}_{2}$, $\psi_{1}$ and
$\overline{\psi}_{2}$ we obtain%
\begin{gather}
e^{-i\rho_{k}x}\left(  z_{k}+1\right)  \sum_{n=0}^{N-1}\left(  -z_{k}\right)
^{n}b_{1,n}(x)-\mathbf{a}(\rho_{k})e^{-i\rho_{k}x}\left(  \overline{z}%
_{k}+1\right)  \sum_{n=0}^{N-1}\left(  -\overline{z}_{k}\right)  ^{n}%
\overline{a}_{2,n}(x)\nonumber\\
-\mathbf{b}(\rho_{k})e^{i\rho_{k}x}\left(  z_{k}+1\right)  \sum_{n=0}%
^{N-1}\left(  -z_{k}\right)  ^{n}a_{1,n}(x)=\left(  \mathbf{a}(\rho
_{k})-1\right)  e^{-i\rho_{k}x},\quad k=1,\ldots,K, \label{sys11}%
\end{gather}%
\begin{gather}
e^{i\rho_{k}x}\left(  \overline{z}_{k}+1\right)  \sum_{n=0}^{N-1}\left(
-\overline{z}_{k}\right)  ^{n}\overline{b}_{2,n}(x)+\widetilde{\mathbf{a}%
}(\rho_{k})e^{i\rho_{k}x}\left(  z_{k}+1\right)  \sum_{n=0}^{N-1}\left(
-z_{k}\right)  ^{n}a_{1,n}(x)\nonumber\\
-\overline{\mathbf{b}}(\rho_{k})e^{-i\rho_{k}x}\left(  \overline{z}%
_{k}+1\right)  \sum_{n=0}^{N-1}\left(  -\overline{z}_{k}\right)  ^{n}%
\overline{a}_{2,n}(x)=\overline{\mathbf{b}}(\rho_{k})e^{-i\rho_{k}x},\quad
k=1,\ldots,K. \label{sys12}%
\end{gather}
Equalities (\ref{eq3}) and (\ref{eq4}) take the form
\begin{equation}
e^{-i\rho_{m}x}\left(  z_{m}+1\right)  \sum_{n=0}^{N-1}\left(  -z_{m}\right)
^{n}b_{1,n}(x)-c_{m}e^{i\rho_{m}x}\left(  z_{m}+1\right)  \sum_{n=0}%
^{N-1}\left(  -z_{m}\right)  ^{n}a_{1,n}(x)=-e^{-i\rho_{m}x}, \label{sys13}%
\end{equation}%
\begin{equation}
e^{i\overline{\rho}_{m}x}\left(  \overline{z}_{m}+1\right)  \sum_{n=0}%
^{N-1}\left(  -\overline{z}_{m}\right)  ^{n}\overline{b}_{2,n}(x)-\overline
{c}_{m}e^{-i\overline{\rho}_{m}x}\left(  \overline{z}_{m}+1\right)  \sum
_{n=0}^{N-1}\left(  -\overline{z}_{m}\right)  ^{n}\overline{a}_{2,n}%
(x)=-\overline{c}_{m}e^{-i\overline{\rho}_{m}x}, \label{sys14}%
\end{equation}
$m=1,\ldots,M$.

For any $x\in\mathbb{R}$, (\ref{sys11})-(\ref{sys14}) represent a system of
$2(K+M)$ linear algebraic equations for the $4N$ unknowns $\left\{
a_{1,n}(x),b_{1,n}(x),\overline{a}_{2,n}(x),\overline{b}_{2,n}(x)\right\}
_{n=0}^{N-1}$, so that $N$ should be chosen such that $2N\leq K+M$. This gives
us an overdetermined system. Solving it aims at finding any of the first SPPS
coefficients, for example, $b_{0}(x)$. As we explained at the beginning of
this section, having computed it at a set of points $x_{l}$ allows us to
recover $q(x)$ from (\ref{q and r}).

\section{Numerical examples\label{Sect Numerics}}

The proposed approach can be implemented directly using an available numerical
computing environment. All the reported computations were performed in Matlab
R2025a on an Intel i7-1360P equipped laptop computer and took from several
seconds to several minutes. For all the integrations involved we used the
Newton--Cottes six point integration rule with $1500$ nodes per unit.
Numerical differentiation which is involved in the last step for recovering
$q(x)$ from the first SPPS coefficients was performed numerically with the aid
of spline interpolation (Matlab routine `\emph{spapi}') and differentiation of
the spline (Matlab routine `\emph{fnder}').

\textbf{Example 1. }Consider
\[
q_{0}(x)=-iA\operatorname*{sech}(x)\exp(-i\gamma A\ln\cosh(x)),\quad
A>0,\,\gamma\in\mathbb{R}.
\]
The analytical expressions of the scattering data for this potential are known
\cite{Tovbis et al 2004}, \cite{Trogdon2021}:
\[
\mathbf{a}(\rho)=\frac{\Gamma\left(  \omega\left(  \rho\right)  \right)
\,\Gamma\left(  \omega\left(  \rho\right)  -\omega_{-}-\omega_{+}\right)
}{\Gamma\left(  \omega\left(  \rho\right)  -\omega_{+}\right)  \,\Gamma\left(
\omega\left(  \rho\right)  -\omega_{-}\right)  },
\]%
\[
\mathbf{b}(\rho)=iA^{-1}2^{-i\gamma A}\frac{\Gamma\left(  \omega\left(
\rho\right)  \right)  \,\Gamma\left(  1-\omega\left(  \rho\right)  +\omega
_{-}+\omega_{+}\right)  }{\Gamma\left(  \omega_{+}\right)  \,\Gamma\left(
\omega_{-}\right)  },
\]%
\[
\omega\left(  \rho\right)  =-i\rho-A\gamma\frac{i}{2}+\frac{1}{2},\quad
\omega_{+}=-iA\left(  T+\frac{\gamma}{2}\right)  ,
\]%
\[
\omega_{-}=iA\left(  T-\frac{\gamma}{2}\right)  ,\quad T=\sqrt{\frac
{\gamma^{2}}{4}-1}.
\]
The zeros of $\mathbf{a}(\rho)$ are given by
\[
\rho_{m}=AT-i(m-\frac{1}{2}),\quad m=1,2,\ldots,M\text{\quad with
}M=\left\lfloor \frac{1}{2}+A\left\vert T\right\vert \right\rfloor .
\]
Also $c_{m}=\mathbf{b}(\rho_{m})$. For the numerical test, following
\cite{Trogdon Olver}, we choose $A=1.0$ and $\gamma=0.1$. To compute the
scattering data, $160$ SPPS coefficients for each Jost solution were
calculated following the recurrent integration procedures devised in Section
\ref{Sect SPPS ZS}. In Fig. 1 the decay of the coefficients $a_{1,n}(0)$,
$a_{2,n}(0)$ is presented.%

\begin{figure}
[ptb]
\begin{center}
\includegraphics[
height=3.4229in,
width=5.3394in
]%
{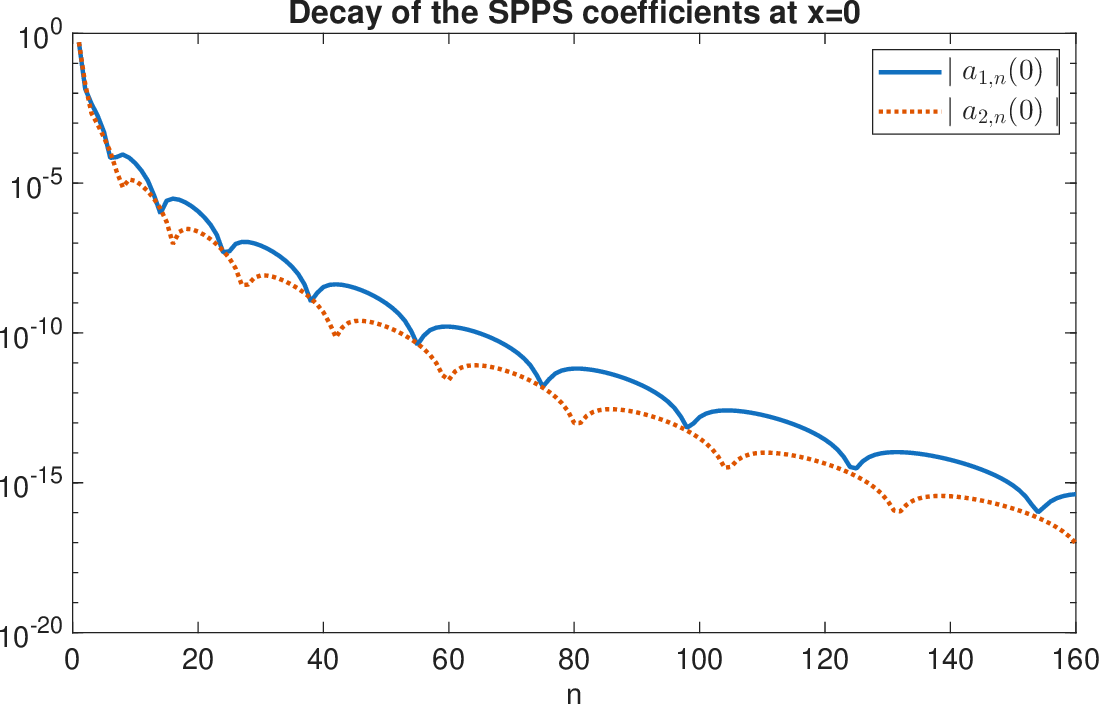}%
\caption{Decay of the computed coefficients $a_{1,n}(0)$ and $a_{2,n}(0)$
corresponding to Example 1, when $n$ increases.}%
\label{Fig1}%
\end{center}
\end{figure}
%

\begin{figure}
[ptb]
\begin{center}
\includegraphics[
height=3.6711in,
width=5.7276in
]%
{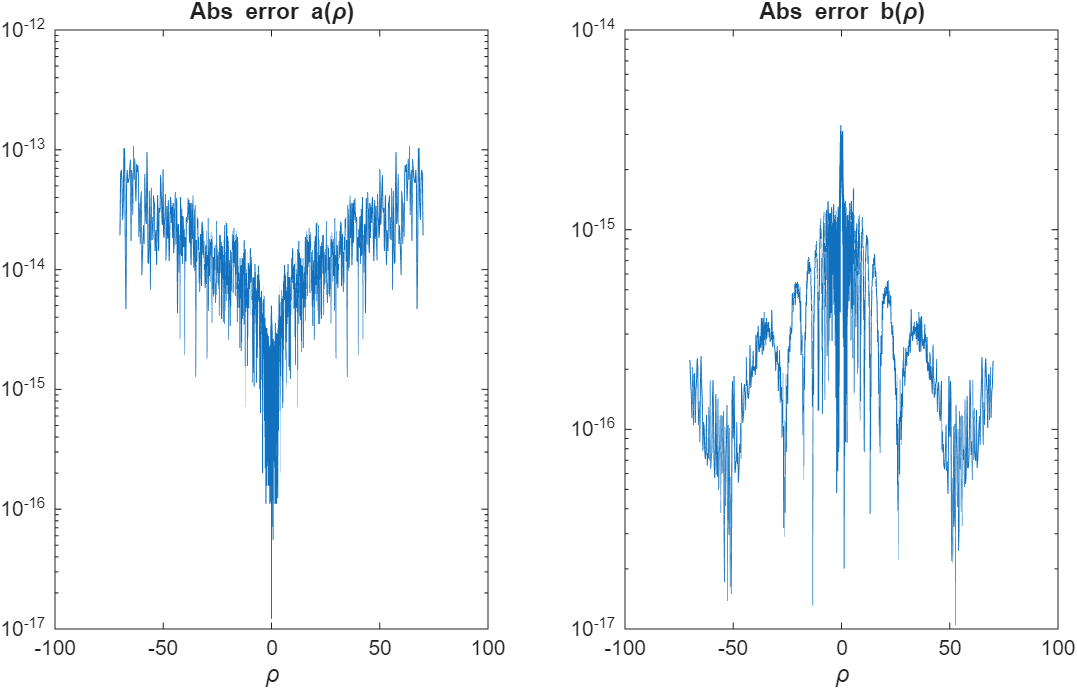}%
\caption{Absolute error of the computed $\mathbf{a}(\rho)$ (left) and
$\mathbf{b}(\rho)$ (right). The maximum absolute error of $\mathbf{a}(\rho)$
resulted in $1.06\cdot10^{-13}$, and that of $\mathbf{b}(\rho)$ in
$3.3\cdot10^{-15}$.}%
\label{Fig2}%
\end{center}
\end{figure}
In Fig. 2 we compare the computed functions $\mathbf{a}(\rho)$ and
$\mathbf{b}(\rho)$ with their analytical expressions. The resulting absolute
error is of order $10^{-13}$. Let us emphasize that the SPPS series
representations allow us to compute the scattering coefficients on immense
intervals of $\rho\in\mathbb{R}$ and in no time, with a non-deteriorating
accuracy and at an arbitrarily large number of points. This is because having
computed once the SPPS coefficients at $x=0$, all further computations reduce
to computing values of the resulting polynomials in $z$ and $\overline{z}$,
when $z$ runs along the unit circle (which corresponds to $\rho$ running along
the real axis).

The values of $\mathbf{a}(\rho_{k})$ and $\mathbf{b}(\rho_{k})$ are computed
at $5000$ points $\rho_{k}$ chosen as follows. Half of them are taken in the
form $\rho_{k}=10^{\alpha_{k}}$ with $\alpha_{k}$ being distributed uniformly
on $\left[  \log(0.001),\log(70)\right]  \approx\left[  -3,\,1.845\right]  $.
This is a logarithmically spaced point distribution on the segment
$[0.001,70]$, with points more densely distributed near $\rho=0.001$ and less
densely distributed as $\rho$ increases. The other half is taken symmetrically
with respect to zero. This logarithmically spaced point distribution for the
computed scattering data may lead to more accurate results when solving the
inverse problem.

One eigenvalue $\rho_{1}\approx0.498749217771909i$ was computed with the
absolute error $2.7\cdot10^{-16}$. We mention that for this, the location of
the roots of the polynomial of the variable $z$ (see Section \ref{SectDirect})
was performed with the aid of the Matlab routine `\emph{roots}'. The
corresponding norming constant $c_{1}\approx
-0.0192926642392854-0.999813879232805i$ was computed with the absolute error
$1.7\cdot10^{-15}$.\ 

Next, the inverse scattering problem is solved for $t=0$, $t=1$ and $t=10$.
The parameter $N$ in (\ref{sys11})-(\ref{sys14}) was chosen as $N=50$. For
$t=0$, we compare the exact potential with the recovered from the computed
scattering data, see the top row in Fig. 3. Here the maximum absolute error of
the recovered potential resulted in $2.1\cdot10^{-4}$. For $t=1$ and $t=10$
the graphs of the computed solution reproduce those from \cite{Trogdon Olver}
(the second and third rows in Fig. 3).%

\begin{figure}
[ptb]
\begin{center}
\includegraphics[
height=4.3638in,
width=6.5743in
]%
{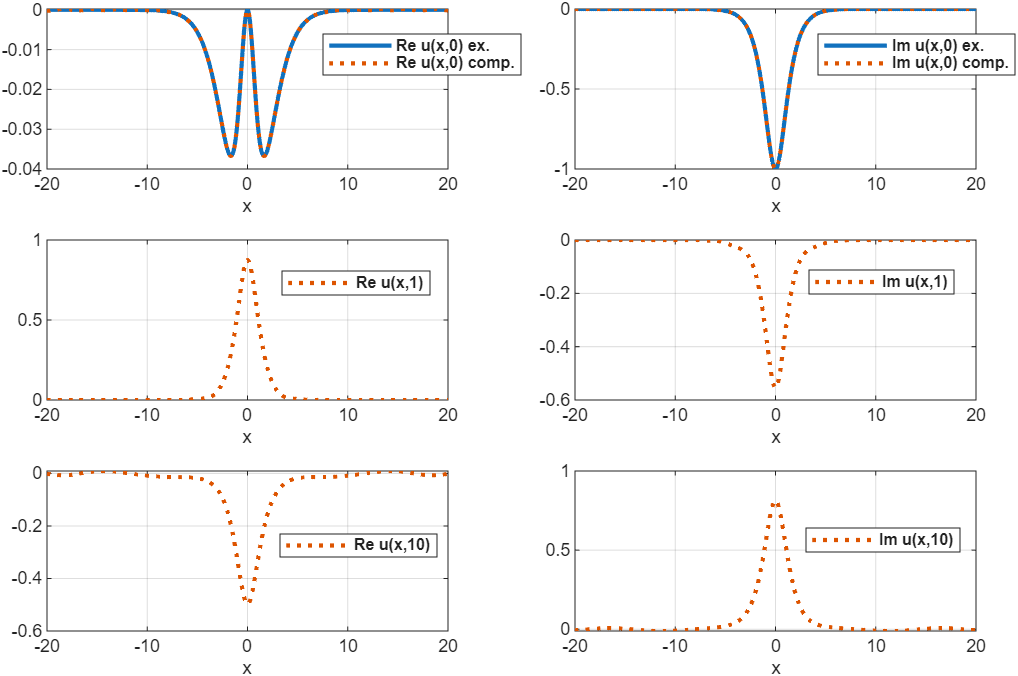}%
\caption{Solution of (\ref{NLSE intro}), (\ref{initial cond}) with $q_{0}(x)$
from Example 1, computed for $t=0$ (top), $t=1$ (middle), $t=10$ (bottom).
Left and right columns show real and imaginary parts, respectively. }%
\label{Fig3}%
\end{center}
\end{figure}

\bigskip

\textbf{Example 2. }A solution of (\ref{NLSE intro}) known in a closed form is
given by \cite[p. 157]{Lamb} (a similar example is considered in \cite{Arico
et al 2011})%
\[
u(x,t)=2\beta\operatorname*{sech}(2\beta x+8\alpha\beta t-\delta
)\exp(-2i\alpha x-4i(\alpha^{2}-\beta^{2})t-i\theta).
\]
This is a single-soliton solution, the constant $2\beta$ defines its
amplitude. The corresponding initial condition (\ref{initial cond}) then
contains
\begin{equation}
q_{0}(x)=2\beta\operatorname*{sech}(2\beta x-\delta)\exp(-i(2\alpha
x+\theta)). \label{q0 Ex2}%
\end{equation}
For the numerical test we chose the values
\[
\alpha=0.5,\quad\beta=\frac{\pi}{2},\quad\delta=0.1,\quad\theta=0.1.
\]
In the first step $60$ SPPS coefficients for each Jost solution were
calculated following the recurrent integration procedures devised in Section
\ref{Sect SPPS ZS}. Here the absolute value of the coefficients decays faster
than in Example 1. The scattering coefficients were computed at points
$\rho_{k}$ distributed as in Example 1. The single eigenvalue $\rho_{1}%
=\alpha+i\beta$ resulted computed with the absolute error $4.6\cdot10^{-14}%
$.\ The value of the corresponding computed norming constant was $c_{1}%
\approx-1.09964966682947-0.110332988730178i$.

We computed the solution of (\ref{NLSE intro}), (\ref{initial cond}) for
$t=0,1,2$. The parameter $N$ in (\ref{sys11})-(\ref{sys14}) was chosen as
$N=50$. In Fig. 4 the absolute value of the solution is shown for the three
instants of time, while in Fig. 5 we compare the exact solution with the
computed one: the real and imaginary parts. It is worth mentioning that the
maximum absolute error of $q(x,0)$ resulted in $2.09\cdot10^{-6}$. That is,
the potential was recovered with such a remarkable accuracy from the computed
scattering data. Even more remarkable is the fact that both for $t=1$ and
$t=2$ the accuracy of the solution did not deteriorate, the maximum absolute
error resulted in $2.13\cdot10^{-6}$ in both cases.%

\begin{figure}
[ptb]
\begin{center}
\includegraphics[
height=3.2422in,
width=5.1967in
]%
{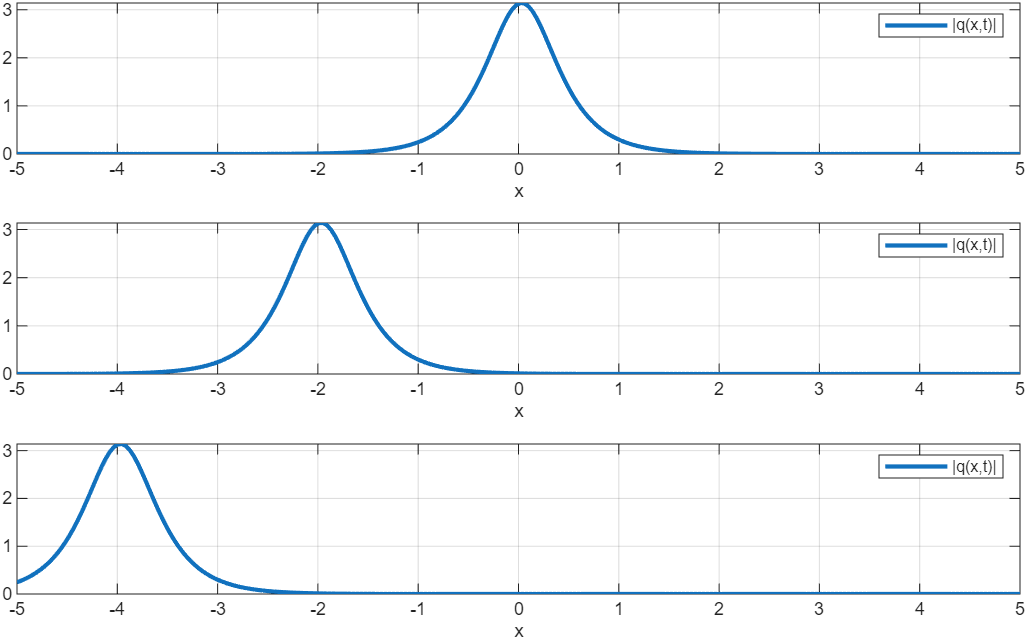}%
\caption{Absolute value of the solution of (\ref{NLSE intro}),
(\ref{initial cond}) with $q_{0}(x)$ defined by (\ref{q0 Ex2}). The solution
is computed for $t=0$ (top), $t=1$ (middle), $t=2$ (bottom).}%
\label{Fig4}%
\end{center}
\end{figure}

\bigskip%

\begin{figure}
[ptb]
\begin{center}
\includegraphics[
height=3.6296in,
width=6.1583in
]%
{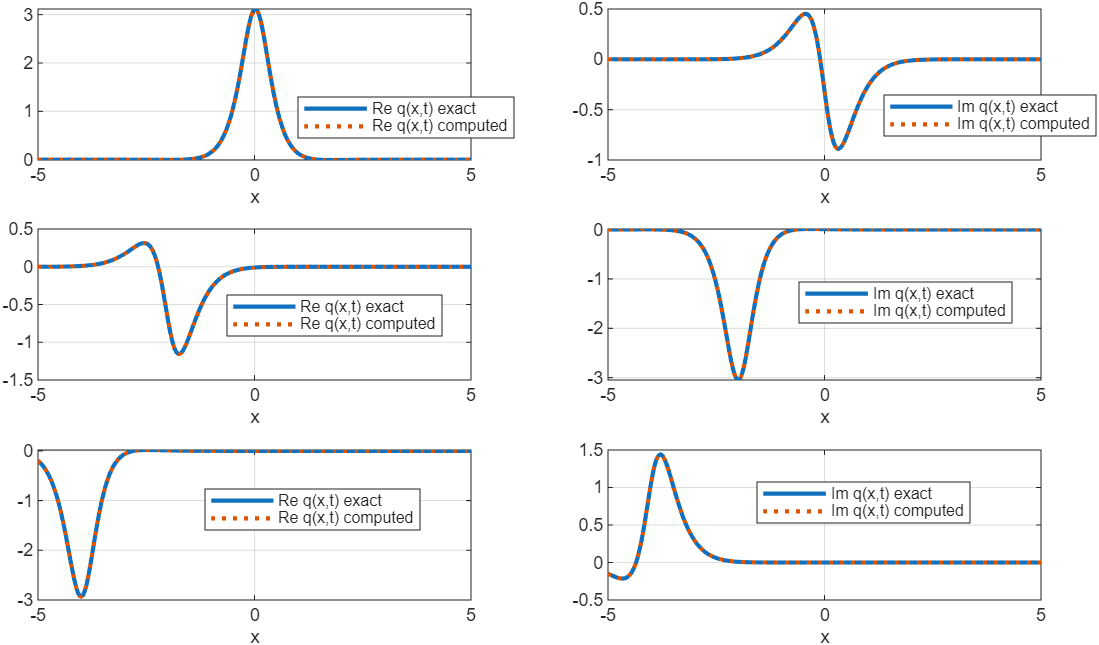}%
\caption{Solution of (\ref{NLSE intro}), (\ref{initial cond}) with $q_{0}(x)$
defined by (\ref{q0 Ex2}), computed for $t=0$ (top), $t=1$ (middle), $t=2$
(bottom). Left and right columns show real and imaginary parts, respectively.
Maximum absolute error of $q(x,0)$, $q(x,1)$ and $q(x,2)$ resulted in
$2.09\cdot10^{-6}$, $2.13\cdot10^{-6}$ and $2.13\cdot10^{-6}$, respectively.}%
\label{Fig5}%
\end{center}
\end{figure}

\textbf{Example 3. }Consider the case \cite{Fermo et al 2015}, \cite{Wahls et
al 2013}
\[
q_{0}(x)=Ae^{i\mu x}e^{-\frac{x^{2}}{\sigma}},
\]
where $A>0$, $\sigma>0$, $\mu\in\mathbb{R}$. Following \cite{Fermo et al 2015}
we choose the values
\[
A=2.5,\quad\sigma=2,\quad\mu=1.
\]
In \cite{Fermo et al 2015} two eigenvalues are reported having their real
parts equal to $-\frac{1}{2}$ and their imaginary parts approximately $1.97$
and $0.79$. Indeed, we obtained%
\begin{align*}
\rho_{1}  &  \approx-0.500000000000079+1.97126262533634i,\\
\rho_{2}  &  \approx-0.499999999999999+0.792849539875588i
\end{align*}
with the corresponding norming constants%
\[
c_{1}\approx-0.999999999999774,\quad c_{2}\approx1.00000000000002.
\]

Solving the inverse scattering problem for the evolved scattering data we
computed $q(x,t)$ for $t=0$ and seven subsequent values of $t$ with the step
$0.3$: $t=0.3,0.6,\ldots,2.1$. The maximum absolute error of the recovered
potential for $t=0$ resulted in $1.3\cdot10^{-3}$. The absolute value of the
solution is presented in Fig. 6.%

\begin{figure}
[ptb]
\begin{center}
\includegraphics[
height=3.9963in,
width=6.6193in
]%
{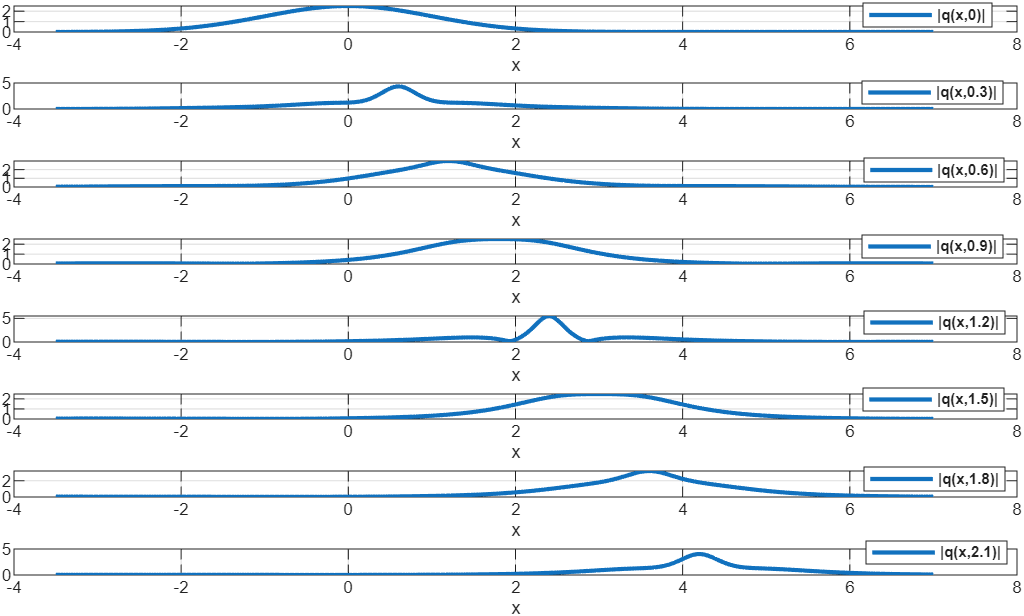}%
\caption{Absolute value of the solution of (\ref{NLSE intro}),
(\ref{initial cond}) for $q_{0}(x)$ from Example 3.}%
\label{Fig6}%
\end{center}
\end{figure}

Some approaches may experience difficulties with less localized initial
conditions. For example, this problem is reported in \cite{Trogdon Olver},
where the authors establish that this happens because the corresponding
reflection coefficient is more oscillatory, which makes it more difficult to
solve the Riemann-Hilbert problem. In the next example we consider such a less
localized initial condition. The method performs well.

\textbf{Example 4. }Consider the case%
\[
q_{0}(x)=\frac{Ae^{i\mu x}}{\left(  x+i\right)  ^{4}}%
\]
with $A=\frac{\pi}{2}$, $\mu=1$. To the difference from the previous examples,
when solving the direct problem and constructing the SPPS coefficients here a
larger interval in $x$ should be considered. In this case the corresponding
interval was chosen as $\left[  -200,200\right]  $ (while in the previous
examples it was sufficient to consider the interval at most $\left[
-12,12\right]  $). This choice gave us a sufficient decay of $\left\vert
q_{0}(x)\right\vert $ at the endpoints of the interval. All other steps did
not require any modification. With $250$ SPPS coefficients computed (in Fig. 7
their decay at $x=0$ is shown) the scattering coefficients were computed at
the points $\rho_{k}$ distributed as in the previous examples. As an
illustration (and checking) of their accuracy, in Fig. 8 we show the accuracy
of the fulfillment of identity (\ref{aatil}).

A single eigevalue was detected
\[
\rho_{1}\approx-2.205978998465+0.485112496978116i
\]
with the corresponding norming constant $c_{1}\approx-1.00000000000001$.

These initial scattering data were used for recovering $q_{0}(x)$ and for
computing the solution $q(x,t)$ for $t=1.2$ and $t=2.5$. The result is
presented in Fig. 9. The maximum absolute error of the recovered $q_{0}(x)$
resulted in $1.9\cdot10^{-2}$.

To control the accuracy of the solution for other values of $t$ we use the
fact that the Wronskian of two solutions of (\ref{ZS1}), (\ref{ZS2}) is
constant. The first SPPS coefficients computed by solving system
(\ref{sys11})-(\ref{sys14}) allow us to compute the Jost solutions
$\varphi(\frac{i}{2},x)$ and $\psi(\frac{i}{2},x)$, see equalities
(\ref{phi(i/2)}). Thus, by checking how $W\left[  \varphi(\frac{i}{2}%
,x);\psi(\frac{i}{2},x)\right]  $ differs from a constant we can estimate the
accuracy of the solution. Thus, for each $t$ we compute the difference
\[
\varepsilon:=\max_{x}\left\vert W\left[  \varphi(\frac{i}{2},x);\psi(\frac
{i}{2},x)\right]  \right\vert -\min_{x}\left\vert W\left[  \varphi(\frac{i}%
{2},x);\psi(\frac{i}{2},x)\right]  \right\vert
\]
as an indicator of the accuracy of the obtained solution. The values of
$\varepsilon$ corresponding to $t=0$, $t=1.2$ and $t=2.5$ resulted
approximately in $0.11$, $0.08$ and $0.03$, thus indicating a satisfactory
level of accuracy.%

\begin{figure}
[ptb]
\begin{center}
\includegraphics[
height=3.1125in,
width=4.8473in
]%
{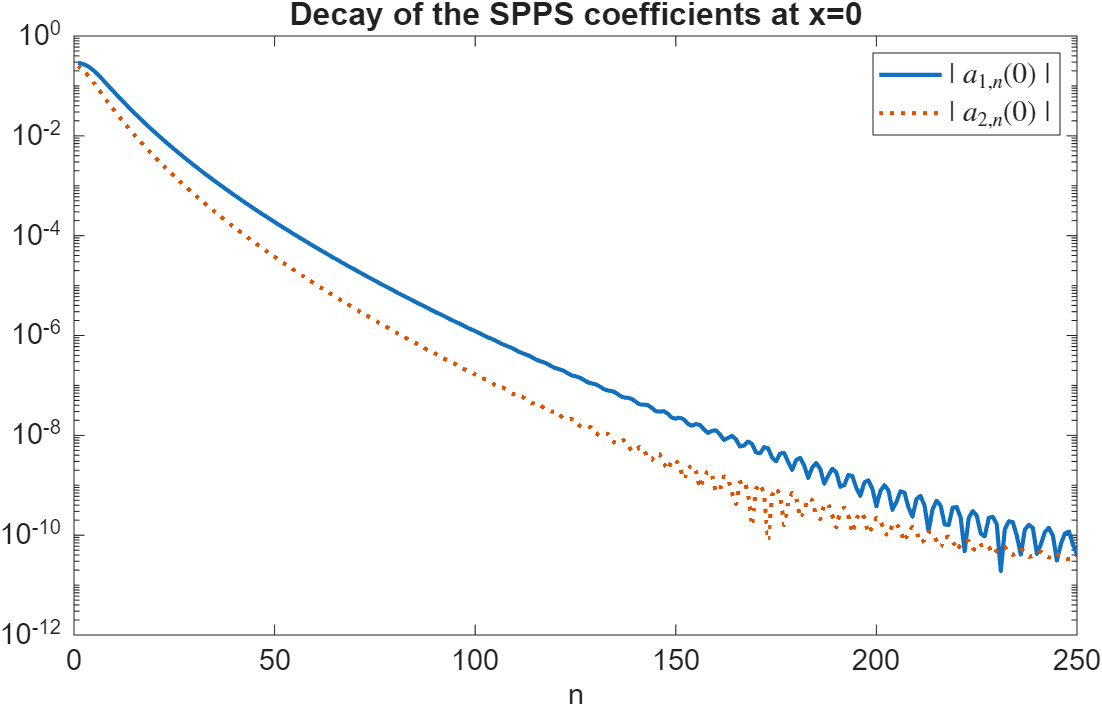}%
\caption{Decay of the computed coefficients $a_{1,n}(0)$ and $a_{2,n}(0)$
corresponding to Example 4, when $n$ increases.}%
\label{Fig7}%
\end{center}
\end{figure}
%

\begin{figure}
[ptb]
\begin{center}
\includegraphics[
height=3.1704in,
width=4.817in
]%
{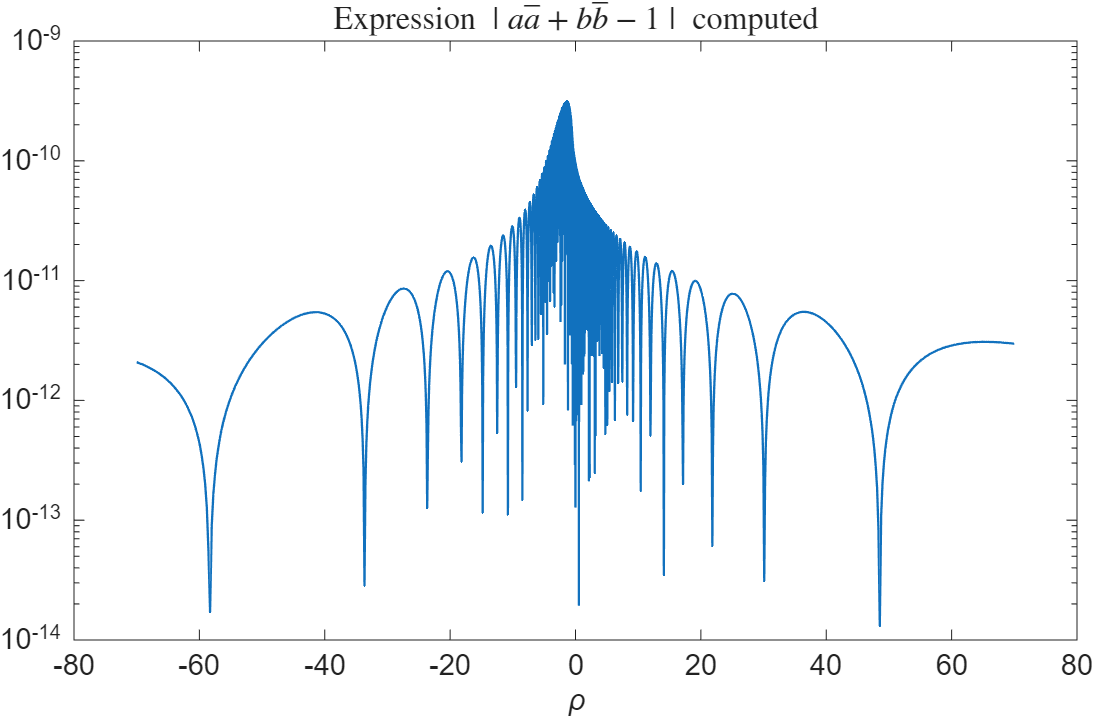}%
\caption{Illustration of the fulfillment of the identity (\ref{aatil}) for the
computed scattering coefficients from Example 4.}%
\label{Fig8}%
\end{center}
\end{figure}

\bigskip%

\begin{figure}
[ptb]
\begin{center}
\includegraphics[
height=3.416in,
width=5.2096in
]%
{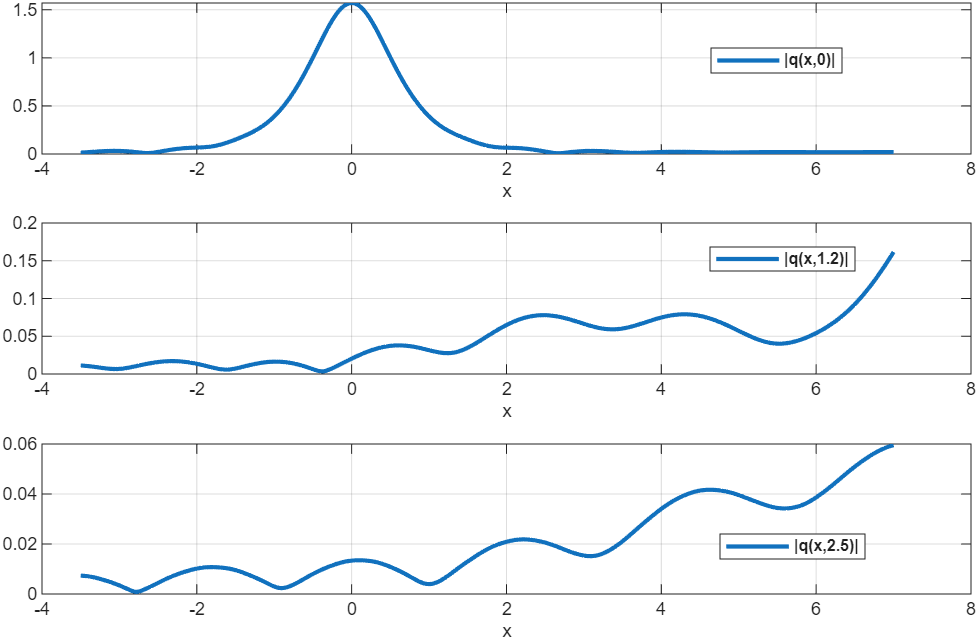}%
\caption{Absolute value of the solution of (\ref{NLSE intro}),
(\ref{initial cond}) with $q_{0}(x)$ from Example 4, computed for $t=0$ (top),
$t=1.2$ (middle) and $t=2.5$ (bottom). }%
\label{Fig9}%
\end{center}
\end{figure}

\section{Conclusions\label{Sect Conclusions}}

A method for solving initial-value problems for the nonlinear Schr\"{o}dinger
equation is developed. It represents an implementation of the inverse
scattering transform method combined with the spectral parameter power series
representations for the Jost solutions for the Zakharov-Shabat system. The
solution of the direct scattering problem reduces to computing the
coefficients of the series, computing products of resulting polynomials and
locating zeros of a polynomial inside the unit disk. In its turn the inverse
scattering problem reduces to the solution of a system of linear algebraic
equations. The solution of the problem for a given instant of time is
recovered from the first elements of the solution vector. The proposed method
leads to an easily implementable, direct and accurate algorithm.

\textbf{Funding information }CONAHCYT, Mexico, grant \textquotedblleft Ciencia
de Frontera\textquotedblright\ FORDECYT - PRONACES/ 61517/ 2020.

\textbf{Data availability} The data that support the findings of this study
are available upon reasonable request.

\textbf{Conflict of interest }This work does not have any conflict of interest.

\end{document}